\documentclass[12pt]{amsart}
\usepackage{amsmath,amstext,amsfonts,amsthm,amssymb,hyperref, amscd,wasysym,stmaryrd}
\usepackage{eucal}
\usepackage[usenames]{color}
\DeclareSymbolFont{bbold}{U}{bbold}{m}{n}
\DeclareSymbolFontAlphabet{\mathbbold}{bbold}
\usepackage[all]{xy}                                        %
\swapnumbers
\newtheorem{thm}[subsubsection]{Theorem}
\newtheorem{lem}[subsubsection]{Lemma}
\newtheorem{prp}[subsubsection]{Proposition}
\newtheorem{crl}[subsubsection]{Corollary}

{  \theoremstyle{definition}

}

\newcommand{\Cat}{\mathtt{Cat}}

\newcommand{\Com}{\mathtt{Com}}

\newcommand{\diag}{\mathrm{diag}\,}

\newcommand{\Env}{\mathtt{Env}}

\newcommand{\Fin}{\mathit{Fin}}
\newcommand{\Fun}{\operatorname{Fun}}
\newcommand{\Funop}{\mathtt{Funop}}

\newcommand{\Hom}{\mathrm{Hom}} 

\newcommand{\Map}{\operatorname{Map}}
\mathchardef\mhyphen="2D

\newcommand{\one}{\mathbbold{1}}

\newcommand{\Op}{\mathtt{Op}}

\newcommand{\Triv}{\mathtt{Triv}}

\newcommand{\Alg}{\mathtt{Alg}}

\newcommand{\cC}{\mathcal{C}}

\newcommand{\cF}{\mathcal{F}}

\newcommand{\cM}{\mathcal{M}}

\newcommand{\cO}{\mathcal{O}}
\newcommand{\cP}{\mathcal{P}}
\newcommand{\cQ}{\mathcal{Q}}
\newcommand{\cS}{\mathcal{S}}

\begin{document}

\title[]{On diagrams of algebras}
\author{Vladimir Hinich}
\address{Department of Mathematics, University of Haifa, Mount Carmel, Haifa 3498838,  Israel}
\email{vhinich@gmail.com}

\begin{abstract}
We present a proof of the formula \cite{L.HA}, 2.4.3.18 for the operad
governing $K$-diagrams of $\cO$-algebras. 
\end{abstract}
 
\maketitle

\section{Introduction}

This is a note about $\infty$-operads. In this note we will use
the word ``category'' to denote $\infty$-categories and ``operad'' to 
denote an $\infty$-operad as defined by Lurie in \cite{L.HA}, Section 2.

To work in a well-defined context, we  accept  quasicategories 
as a model for $\infty$-categories; but
all our constructions are presented in a $\infty$-categorical language, 
as it is described in \cite{H.EY}, Section 2, so that they make sense in 
any model.

The term ``conventional categories'' stands for those  categories whose
spaces of morphisms are equivalent to sets.

In this note $\Cat$ denotes the category of small categories, $\Fin_*$
is the category of finite pointed sets and an operad is a functor 
$p:\cO\to\Fin_*$ satisfying the standard properties of Definition 2.1.1.10 of~\cite{L.HA}. In particular, $\Com=\Fin_*$ is the operad for commutative 
algebras.

The category of operads $\Op$ is defined as the subcategory of $\Cat_{/\Fin_*}$ spanned by the operads, with the arrows preserving cocartesian liftings of the inerts. It can also be defined as a
 Bousfield localization as follows. 
Let $\Cat^+_{/\Fin_*^\natural}$ the category  of marked categories
over $\Fin_*$ endowed with the standard marking (inert arrows are marked).
Then $\Op$ indentifies with the full subcategory of 
$\Cat^+_{/\Fin_*^\natural}$ spanned by the operads with the inerts
as the marked arrows. The full embedding 
$R:\Op\to\Cat^+_{/\Fin_*^\natural}$ admits a left adjoint
$$
L:\Cat^+_{/\Fin_*^\natural}\to\Op,
$$
so that $\Op$ becomes the Bousfield localization of 
$\Cat^+_{/\Fin_*^\natural}$ with respect to the equivalence
determined by $L$ (called the operadic equivalence).

\subsection*{Acknowledgement}
The author is very grateful to I.~Moerdijk who pointed out to a gap in the 
original proof. The operad  $\cC^K$ defined in this note, was mentioned in 4.1.3 of \cite{HM}.
 
\section{An operad for diagrams of algebras}

Let $\cO$ and $\cC$ be operads and $K$ be a category. In this note we
discuss the functor
assigning to $\cO,\cC$ and $K$  the category $\Fun(K,\Alg_\cO(\cC))$.
This functor is representable in different ways.
\begin{itemize}
\item As a functor of $\cC$, it is(co)represented by an operad that we denote by $\cO_K$, see~\cite{H.EY}, 2.10.5(3). This is the operad governing $K$-diagrams of $\cO$-algebras. By definition, $\cO_K$ is an operad endowed with an operadic equivalence
$$
\gamma:K\times\cO\to\cO_K,
$$
where $K\times\cO$ is considered as marked category over 
$\Fin_*^\natural$, where an arrow $(\alpha,\beta)$ in $K\times\cO$ is marked iff $\alpha$ is an equivalence and $\beta$ is inert. 
\item As a functor of $\cO$, it is represented by the operad $\cC^K$
so that, in the case when $\cC$ is a symmetric monoidal category,
$\cC^K$ is the symmetric monoidal category of functors $K\to\cC$, see
Section~\ref{sec:path}. 
\end{itemize}

Furthermore, both $\cO_K$ and $\cC^K$ have an explicit expression in 
terms of the operad $K^\sqcup$ defined in \cite{L.HA}, 2.4.3.

Here are the main results of this work.
\begin{itemize}
\item[1.] The operad $K^\sqcup$ is flat, see Lemma~\ref{lem:ksqcup-flat}.
\item[2.] There is an equivalence $\cC^K=\Funop(K^\sqcup,\cC)$,
see~\ref{crl:path}.
\item[3.] There is an equivalence $\cO_K=\cO\times_{\Com}K^\sqcup$, see
Theorem~\ref{thm:OE} proven in~\ref{ss:proof-OE}.
\end{itemize}

The most interesting equivalence is the Claim 3. It was first mentioned in~\cite{L.HA}, 2.4.3.18, but the reasoning there was based on an incorrect Remark 2.4.3.6. 

\subsection{}
Recall the definition of $K^\sqcup$, \cite{L.HA}, 2.4.3.1.

Define $\Gamma^*$ as the  (conventional) category of pairs 
$(I_*,i)$ with $I_*\in\Fin_*$ and $i\in I$, with the arrows
$(I_*,i)\to(J_*,j)$ given by arrows $I_*\to J_*$ carrying $i$ to $j$.
The functor $\pi:\Gamma^*\to\Fin_*$ carries
$(I_*,i)$ to $I_*$.

For $K\in\Cat$, we define $K^\sqcup$ as a category over 
$\Com=\Fin_*$ representing the functor
$$
B\mapsto\Map(B\times_{\Fin_*}\Gamma^*,K).
$$

The fiber of $K^\sqcup$ at $I_*\in\Com$ is $K^I$; an arrow
in $K^\sqcup$ over $\alpha:I_*\to J_*$ from $x:I\to K$ to $y:J\to K$
is given by a collection of arrows $x(i)\to y(j)$ for all pairs $(i,j)\in
I\times J$ with $\alpha(i)=j$.

\subsubsection{} 
In the case when $K$ is a conventional category, $K^\sqcup$ is a conventional operad. Its colors are the objects of $K$ and 
an operation from $\{x_i\}$ to $y$ is given by a collection of arrows
$x_i\to y$. The composition of operations is defined in an obvious 
way.

\subsubsection{}
In the special case when $K\in\Cat$ has finite coproducts, $K^\sqcup$
is the operadic presentation of the cocartesian SM category $K$,
see~\cite{L.HA}, 2.4.3.12.

\subsection{Flatness of $K^\sqcup$}

Recall \cite{H.EY}, 2.8.2, that an operad $\cO$ is called {\sl flat}
if for any pair of composable active arrows $s:x_0\to x_1\to x_2$ in $\Fin_*$ the base change $\cO\times_{\Fin_*}[2]\to[2]$ is flat
in the sense of~\cite{L.HA}, B.3. If an operad $\cO$ is flat, one can
define a functor $\cP\mapsto\Funop(\cO,\cP)$ so that
$$
\Alg_\cQ(\Funop(\cO,\cP))=\Alg_{\cQ\times\cP}(\cP),
$$
where $\cQ\times\cO$ is the product in $\Op$. 

\begin{lem}
\label{lem:ksqcup-flat}
The operad $K^\sqcup$ is flat.
\end{lem}
\begin{proof}
By \cite{H.EY}, 2.8.2, we have to verify that for any pair of composable active arrows $s:x_0\to x_1\to x_2$ in $\Fin_*$ the base change $K^\sqcup\times_{\Fin_*}[2]\to[2]$ is flat.

 Since the restriction of
$K^\sqcup$ to the active part of $\Fin_*$ is a cartesian fibration,
the flatness is immediate.
\end{proof}

The following result slightly generalizes~\cite{H.EY}, 2.8.10.

\begin{prp}
\label{prp:Kalgebras}
Let $\cM$ be a  symmetric monoidal category and $\cM^\otimes$
be its operadic presentation. Then for any $K\in\Cat$ the operad
$\Funop(K^\sqcup,\cM^\otimes)$ is the operadic presentation of the
symmetric monoidal category $\Fun(K,\cM)$. Moreover, if $\cM$ is cartesian, $\Funop(K^\sqcup,\cM^\otimes)$ is also cartesian.
\end{prp}
\begin{proof}
We denote $\cF=\Funop(K^\sqcup,\cM^\otimes)\in\Op$. Let us first of all describe the underlying category $\cF_1$. One has
$$
\cF_1=\Alg_{\Triv}(\cF)=\Alg_{\Triv\times_{\Fin_*}K^\sqcup}(\cM^\otimes)=\Fun(K,\cM).
$$
Here $\Triv$, the trivial operad, is the subcategory of $\Fin_*$ spanned by the inert arrows.

Let $f=(f_1,\ldots,f_n)$ and $g$ be functors $K\to\cM$. Let us describe
$\Map^p(f,g)$, the space of arrows in $\cF$ over the active arrow 
$p:\langle n\rangle\to\langle 1\rangle$. The calculation is very similar to (but considerably easier) \cite{H.EY}, 4.2.2 and 4.2.3.
We denote by $C_n$ the operad generated by one $n$-ary operation
and by $C_n^\circ$ the subcategory of inert arrows in $C_n$. The
pair $(f,g)=(f_1,\ldots,f_n,g)$ defines a $C_n^\circ$-algebra in $\cF$
and the space $\Map^p_\cF(f,g)$ is  the fiber of
the restriction map
$$
\Alg_{C_n}(\cF)\to\Alg_{C_n^\circ}(\cF)
$$
at $(f,g)$.
Note that $\Alg_{C_n}(\cF)=\Alg_{C_n\times_{Fin_*}K^\sqcup}(\cM^\otimes)$.

Denote by $K^\sqcup_p$, $\cM^\otimes_p$ the categories over $[1]$ obtained
from $K^\sqcup$, $\cM^\otimes$ by the base change $[1]\to\Fin_*$ defined
by the active arrow $p:\langle n\rangle\to\langle 1\rangle$. One has
$$
\Alg_{C_n\times_{Fin_*}K^\sqcup}(\cM^\otimes)=\Fun_{[1]}(K^\sqcup_p,
\cM^\otimes_p).
$$
Now, $K^\sqcup_p$ is the cartesian fibration classified by the
diagonal map $K\to K^n$, whereas $\cM^\otimes_p$ is the cocartesian
fibration classified by the (multiple) tensor product $\cM^n\to\cM$.
This allows one to identify $\Map^p(f,g)$ with the space
$\Hom_{\Fun(K,\cM)}(f_1\otimes\ldots\otimes f_n,g)$
where $f_1\otimes\ldots\otimes f_n$ is defined as the composition
$$
K\stackrel{\diag}{\to} K^n\stackrel{\prod f_i}{\longrightarrow}\cM^n\to\cM,
$$
where the last map is the (multiple) tensor product in $\cM$.
This proves that $\Funop(K^\sqcup,\cM^\otimes)$ is a symmetric monoidal category.  

Let now $\cM^\otimes$ be cartesian. Let $1\in\cM$
and $\one\in\Fun(K,\cM)$ be final objects of $\cM$ and of $\Fun(K,\cM)$ 
respectively.
Given $f_1,f_2:K\to\cM$, we have to verify that the diagram
$$f\otimes\one\longleftarrow f\otimes g\longrightarrow \one\otimes g$$
is cartesian that is equivalent to saying that its evaluation at any $x\in K$ is cartesian in $\cM$. This follows from the fact that $\cM$ is cartesian.

\end{proof}

\subsection{}

The natural map $\gamma:K\times\Com\to K^\sqcup$ is given by the 
projection
$$
(K\times\Com)\times_{\Com}\Gamma^*=K\times\Gamma^*\to K.
$$

Given an operad $\cO$, we obtain, by the base change, the map
\begin{equation}
\gamma_\cO:K\times\cO\to K^\sqcup\times_{\Com}\cO.
\end{equation}

We see $K\times\cO$ as an object of the category 
$\Cat^+_{/\Fin^\natural_*}$. 

The following result is central for our discussion.

\begin{thm}
\label{thm:OE}
$\gamma_\cO$ is an operadic equivalence.
\end{thm}

The theorem is proven in~\ref{ss:proof-OE}.

\section{Path space of an operad}
\label{sec:path}

\subsection{}
Given an operad $\cP$ and a category $K$, we define a category 
$\cP^K$ over $\Fin_*$ by the formula
$$
\cP^K=\Fun(K,\cP)\times_{\Fun(K,\Fin_*)}\Fin_*,
$$
where $\Fin_*\to\Fun(K,\Fin_*)$ assigns to $I_*$ the constant functor
$K\to\Fin_*$ with the value $I_*$.
\begin{lem}
$\cP^K$ is an operad.
\end{lem}
\begin{proof}
Let $\cM=\Env(\cP)$ be the symmetric monoidal envelope of $\cP$
and let $\cM_1$ be the underlying category. In this case $\cM^K$ is
a cocartesian fibration over $\Fin_*$ representing the standard symmetric 
monoidal structure on $\Fun(K,\cM_1)$. Obviously $\cP^K$ is the full suboperad of $\cM^K$ spanned by $\cP_1\subset\cM_1$.
\end{proof}

The following result is almost immediate.
\begin{prp}
There is a canonical equivalence
$$
\Map_\Op(\cO_K,\cP)=\Map_\Op(\cO,\cP^K).
$$
\end{prp}
\begin{proof}
One has
\begin{eqnarray}
\nonumber 
\Map_\Op(\cO_K,\cP)=\Map_{\Cat^+_{/\Fin_*^\natural}}(\cO\times
K^\flat,\cP)=\hspace{2in}\\
\nonumber\Map_{\Cat^+_{/\Fin_*^\natural}}(\cO,\Fun(K,\cP)
\times_{\Fun(K,\Fin_*)}\Fin_*)=\Map_\Op(\cO,\cP^K).
\end{eqnarray}
\end{proof}

\subsection{Proof of (\ref{thm:OE})}
\label{ss:proof-OE}
The map $\gamma_\cO$ induces a map of operads
$$
\bar\gamma_\cO:\cO_K\to\cO\times_\Com K^\sqcup.
$$
We will prove it is an equivalence of operads using the reconstruction
theorem \cite{HM}, 4.4.4.

The map $\bar\gamma_\cO$ is equivalence on colors, so it is enough to verify
that it induces the equivalence of the categories of algebras in $\cS$.
By Proposition~\ref{prp:Kalgebras}
$$
\Alg_\cO(\Funop(K^\sqcup,\cS))=\Alg_\cO(\Fun(K,\cS))=
\Fun(K,\Alg_\cO(\cS)).
$$
The map $\gamma_\cO$ identifies this category with 
$\Alg_{K^\flat\times\cO}(\cS)$. This proves the result.

Finally, one has
\begin{crl}
\label{crl:path}
There is a natural equivalence $\cC^K=\Funop(K^\sqcup,\cC)$.
\end{crl}
\begin{proof}
By Theorem~\ref{thm:OE} the operadic equivalence $\gamma_\cO$ induces
an equivalence 
$$
\Alg_{\cO\times_\Com K^\sqcup}(\cC)\to\Alg_{\cO_K}(\cC)=\Alg_\cO(\cC^K),
$$
or, in other words, an equivalence
$$
\Alg_{\cO}(\Funop(K^\sqcup,\cC))\to\Alg_\cO(\cC^K).
$$
This implies the claim.
\end{proof}

\end{document}